\documentclass{endm}
\usepackage{endmmacro}
\usepackage{graphicx, amsmath}

\newtheorem{question}{Question}

\newcommand\cG{\mathcal G}
\newcommand\cM{\mathcal M}
\newcommand\Forb{\operatorname{Forb}}

\begin{document}

\begin{verbatim}\end{verbatim}\vspace{2.5cm}

\begin{frontmatter}

\title{Forbidden Induced Subgraphs}

\author{Thomas Zaslavsky\thanksref{myemail}}
\address{Department of Mathematical Sciences\\ Binghamton University (State University of New York)\\ Binghamton, NY 13902-6000, U.S.A.}
\thanks[myemail]{Email:
   \href{mailto:zaslav@math.binghamton.edu} {\texttt{\normalshape
   zaslav@math.binghamton.edu}}} 

\begin{abstract}
In descending generality I survey: five partial orderings of graphs, the induced-subgraph ordering, and examples like perfect, threshold, and mock threshold graphs.  The emphasis is on how the induced subgraph ordering differs from other popular orderings and leads to different basic questions.
\end{abstract}

\begin{keyword}
Partial ordering of graphs, hereditary class, induced subgraph ordering, perfect graph, mock threshold graph.
\end{keyword}

\end{frontmatter}

\section{Preparation}\label{intro}

This is a very small survey of partial orderings of graphs, hereditary graph classes (mainly in terms of induced subgraphs), and characterizations by forbidden containments and by vertex orderings, leading up to a new graph class and a new theorem.

Our graphs, written $G = (V,E)$, $G' = (V',E')$, etc., will be finite, simple (with some exceptions), and unlabelled; thus we are talking about isomorphism types (isomorphism classes) rather than individual labelled graphs.  Graphs can be partially ordered in many ways, of which five make an interesting comparison.

We list some notation of which we make frequent use:
\begin{quote}
\qquad$N(G;v) =$ the neighborhood of $v$ in $G$.

\quad\;$d(G; v) =$ the degree of vertex $v$ in $G$.

\ \qquad$G{:}X =$ the subgraph of $G$ induced by $X \subseteq V$.
\end{quote}

\section{Five Kinds of Containment}

We say a graph $G$ \emph{contains} $H$ if $H \leq G$, where $\leq$ denotes some partial ordering on graphs.  There are many kinds of containment; each yields a different characterization of each interesting hereditary graph class.

Two kinds of graph that make good models are planar and outerplanar graphs.  A \emph{planar graph} can be drawn in the plane with no crossing edges.  An \emph{outerplanar graph} can be drawn with no crossing edges and with all vertices on the boundary of the infinite region.  These two types can be characterized by forbidden contained graphs.
\medskip

1.  The \emph{subgraph} ordering:  $G' \subseteq G$ means $V' \subseteq V$ and $E' \subseteq E$.  (Technically, this definition applies to labelled graphs.  We apply it to unlabelled graphs by forgetting the labels.).  Examples:
\begin{enumerate}
\item Planarity.  There are infinitely many forbidden subgraphs, i.e., graphs such that $G$ is planar iff it contains none of those graphs.
\item Outerplanarity.  There are also infinitely many forbidden subgraphs.  
\end{enumerate}
But both infinities are repaired by the next form of containment.  \emph{Subdividing} a graph means replacing each edges by a path of positive length (thus, $G$ is a subdivision of itself).
\medskip

2.  The \emph{subdivided subgraph} or \emph{topological subgraph} ordering:  $G' \subseteq_t G$ means there is a subdivision $G''$ of of $G'$ such that $G'' \subseteq G$.  Examples:
\begin{enumerate}
\item Planarity.  By Kuratowski's Theorem there are two forbidden subdivided subgraphs ($K_5$ and $K_{3,3}$, if anyone forgot).
\item Outerplanarity.  There are also two forbidden subdivided subgraphs, $K_4$ and $K_{2,3}$, by Chartrand and Harary \cite{ChaH}.
\end{enumerate}
\medskip

3.  The \emph{minor} ordering:  $G' \preceq G$ means $G$ has a subgraph that contracts to $G'$.  (Contraction means shrinking an edge to a point, thereby combining the endpoints into a single vertex.  Or, several edges may be contracted.)
Examples:
\begin{enumerate}
\item Planarity.  By Wagner's extension of Kuratowski's Theorem there are two forbidden minors ($K_5$ and $K_{3,3}$ again).
\item Outerplanarity.  There are two forbidden minors ($K_4$ and $K_{2,3}$; trivially since any $K_4$ or $K_{2,3}$ minor is a subdivision).
\end{enumerate}
\medskip

4.  The \emph{induced subgraph} ordering:  $G' \sqsubseteq G$ means $V' \subseteq V$ and $E'$ contains all edges of $G$ that have both endpoints in $V'$.  The induced subgraph ordering is a very active research topic currently.  Examples (both obvious):
\begin{enumerate}
\item Planarity.  There are infinitely many forbidden induced subgraphs.
\item Outerplanarity.  Also infinitely many forbidden induced subgraphs.
\end{enumerate}
Does subdivision eliminate the infinities?  We subdivide:
\medskip

5.  The \emph{subdivided induced subgraph} or \emph{topological induced subgraph} ordering:  $G' \sqsubseteq_{t} G$ means there is a subdivision $G''$ of $G'$ such that $G'' \sqsubseteq G$.  This ordering seems to be newly regarded.  Examples:
\begin{enumerate}
\item Planarity.  Infinitely many forbidden topological induced subgraphs.
\item Outerplanarity.  We have not yet decided finiteness of the number of forbidden topological induced subgraphs.
\end{enumerate}

Despite that last uncertainty, subdivision can matter; see Theorem \ref{mtgp}.

\section{Hereditary Classes and Forbidden Subgraphs}

A \emph{hereditary class} of graphs in a partial ordering of graphs, $\leq$, is a class $\cG$ of graphs for which
$G \in \cG\ \&\ G' \leq G \implies G' \in \cG.$
The \emph{forbidden graphs} for $\cG$ are the minimal non-members of $\cG$.  We write $\Forb(\cG)$ for the set of minimal non-members.  Thus, 
$$G \in \cG \iff \text{ no element of } \Forb(\cG) \text{ is } \leq G.$$
The main question:  What is $\Forb(\cG)$?  Is there any general statement about it?  The most famous answer of this kind at present is the Robertson--Seymour Graph Minors Theorem:

\begin{theorem}
In the minor ordering $\preceq$, every $\Forb(\cG)$ is finite.
\end{theorem}

None of the other orderings has this finiteness property.

The Robertson--Seymour graph minors theory is comprehensive, including a structure theory of graphs.  Another aspect of it is recognition.  
\emph{Class P} is the class of problems for which the time required is not more than a polynomial function of the size of the input (these are often called ``fast'' computations).  
\emph{Class NP} is the class of problems for which the required time is polynomial if there is a hint.  Trivially, P (where no hint is needed) $\subseteq$ NP.  The notorious open question:  Is P $\subsetneq$ NP?

\begin{corollary}
Given a fixed graph $H$, the question  ``For every graph $G$, decide whether $H \preceq G$''  is in Class P.
\end{corollary}

This leads to our main problem:

\begin{question}\label{finiteness}
Is anything like these results true in other orderings?  In particular, is there any kind of finiteness for interesting hereditary classes?
\end{question}

See the introduction to \cite{KRT} for a review of questions like this for several orderings, especially the immersion ordering (which I have omitted).

\begin{question}\label{containment}
How difficult is it to decide whether $H \leq G$ in other orderings?
\end{question}

I am especially interested in the \emph{induced} and \emph{topological induced} subgraph orderings.\footnote{For the subgraph and induced-subgraph orderings, Question \ref{containment} is easily answered in polynomial time.}  
They are the subject of the rest of this report.

\section{Holes and Antiholes}

A \emph{hole} in a graph is an induced cycle $C_l$ for $l\geq4$.  An \emph{antihole} is an induced subgraph isomorphic to the complement $\overline C_l$ for $l\geq 4$.
Forbidding holes and antiholes is important.

\subsection{Holes and Triangles:  Forests}

Perhaps someone may be surprised by the following form of statement, an induced-subgraph characterization.

\begin{proposition}
$G$ has no holes and no induced triangles 
$\iff$
$G$ contains no cycles 
$\iff$
$G$ is a forest.
\end{proposition}

For a constructive characterization of forests we have a slightly unconventional restatement of the fact that every tree other than $K_1$ has a leaf:

\begin{proposition}
$G$ is a forest $\iff$ there is a vertex ordering $(v_1,v_2,\ldots,v_n)$ such that
$d(G{:}\{v_1,v_2,\ldots,v_i\}; v_i) \leq 1.$
\end{proposition}

The vertex ordering implies that recognizing forests is in Class P.

\subsection{Holes:  Chordal Graphs (Triangulated Graphs)}

A graph is \emph{chordal} if every cycle in $G$ of length $\geq4$ has a chord.
Restated:

\begin{proposition}
$G$ is chordal $\iff$ $G$ has no holes.
\end{proposition}

There is an important constructive characterization by Dirac \cite{D}:

\begin{proposition}
$G$ is chordal $\iff$ there is a vertex ordering $(v_1,v_2,\ldots,v_n)$ such that
$N(G{:}\{v_1,v_2,\ldots,v_i\}; v_i)\ \text{ is a clique}.$
\end{proposition}

This vertex ordering implies that recognizing chordal graphs is in Class P.

\subsection{Odd Holes:  Bipartite Graphs }

Parity of holes matters.

\begin{proposition}
$G$ has no odd holes and no triangles $\iff$ $G$ is bipartite.
\end{proposition}

Is there a constructive characterization involving vertex ordering?  Apparently not.  Still, bipartite graphs matter:
\begin{enumerate}
\item Many computational problems have simple algorithms for bipartite graphs but difficult ones for non-bipartite graphs.  Many have algorithms in P for bipartite graphs, but not for all graphs.
\item Bipartite graphs are the answer to some good graph-theory questions.
\item Bipartite graphs are good subjects for integer linear programming, but general graphs are not.
\end{enumerate}

Is there an ``anti'' version of this?  
$G$ has no odd antiholes and no independent vertex triples $\iff$ it consists of two cliques and connecting edges.  
So odd antiholes are not very interesting?  Well \dots

\subsection{Odd Holes and Antiholes: Perfect Graphs}

Some more notation:
\begin{quote}
\qquad$\chi =$ chromatic number.

\quad$\omega =$ clique number (the largest size of a clique).
\end{quote}
In general, obviously, $\chi(G) \geq \omega(G)$.  So, when are they equal?  A \emph{perfect graph} is a graph $G$ that satisfies $\chi(G') = \omega(G')$ for every $G' \sqsubseteq G$.  Note that the requirement extends to induced subgraphs.  Why do we care about these graphs?
\begin{enumerate}
\item Perfect graphs are the answer to some good graph-theory questions.
\item Perfect graphs are good for polyhedral combinatorics.
\item Many computational problems have fast algorithms (Class P) for perfect graphs but not for all graphs.
\end{enumerate}

\begin{conjecture}[Berge, 1961]
$G$ is perfect $\iff$ $\overline G$ is perfect $\iff$ $G$ has no odd holes or antiholes.
\end{conjecture}

The Strong Perfect Graph Theorem of Chudnovsky--Robertson--Seymour--Thomas \cite{CRST} says Berge was right.

\begin{theorem}
$G$ has no odd holes and no odd antiholes $\iff$ it is perfect.
\end{theorem}

Graphs with no odd holes or antiholes have become known as \emph{Berge graphs}.  There is now a thorough structural decomposition that enables polynomial-time recognition of perfect graphs (see Chudnovsky \emph{et al.}\ \cite{CS}).  
One direction of current research is to modify the Berge exclusions, hoping that valuable properties of perfect graphs remain valid.  For instance \ldots

\subsection{Even Holes}

A \emph{nearly Berge} graph has no even holes.
One reason for this definition is that having no even holes $\implies$ no induced $C_4$ $\implies$ no antiholes $\overline C_l$ for $l\geq6$.  So, about half the Berge exclusions are preserved.
There is an attractive constructive description by vertex ordering \cite{SV}:

\begin{theorem}
A graph without even holes has a vertex ordering such that $N(v_i)$ in $G{:}\{v_1,v_2,\ldots,v_i\}$ is chordal.
\end{theorem}

But not the converse, unfortunately, so this is not a characterization (see \cite{V}) and it does not give a recognition algorithm.

\section{Three Small Exclusions:  Threshold Graphs}

\emph{Threshold graphs}, introduced by Chv\'atal and Hammer \cite{CH,CH2}, can be defined by forbidden induced subgraphs: they are the graphs such that $G \not\sqsupseteq P_4, C_4, \overline C_4$.  (That is not the original definition, which is more complicated and off topic.)  This characterization implies that the complement of a threshold graph is again a threshold graph.

A constructive characterization is also due to Chv\'atal and Hammer.  (A vertex is \emph{dominating} if it is adjacent to all other vertices.)

\begin{proposition}
$G$ is threshold $\iff$
$\exists$ a vertex ordering such that each $v_i$ is isolated or dominating in $G{:}\{v_1,v_2,\ldots,v_i\}$.
\end{proposition}

The vertex ordering, once again, implies that recognizing threshold graphs is in Class P (and fast).

Is there a generalization?  Yes \dots

\section{Relaxed Vertex Ordering:  Mock Threshold Graphs}

So far we have gone from forbidden induced subgraphs to a vertex ordering.  Now it's reversal time.  We relax the threshold ordering requirement to define a new class of graphs.  In 2014 Sivaraman said a graph is \emph{mock threshold}\footnote{I acknowledge responsibility for this name.} if there is a vertex ordering such that $d(G{:}\{v_1,v_2,\ldots,v_i\};v_i) \leq 1$ or $\geq i-2$.
Two essential properties are:
\begin{enumerate}
\item The complement of a mock threshold graph is mock threshold.
\item An induced subgraph of a mock threshold graph is mock threshold.
\end{enumerate}
Hence the obvious question.

\begin{question}
What are the forbidden induced subgraphs for the class $\cM$ of mock threshold graphs?
\end{question}

Easy forbidden induced subgraphs are the holes $C_5, C_6, C_7, \ldots$ and their antiholes.
Can we find the others?  Are there infinitely many?

\begin{theorem}[Behr--Sivaraman--Zaslavsky--computer \cite{BSZ}]
$\Forb(\cM)$ \linebreak consists of all holes $C_l$ and antiholes $\overline C_l$ for $l\geq5$, and $318$ other graphs.\footnote{We thank Jeffrey Nye for computational help.}
\end{theorem}

Once more, because of the definition by vertex ordering, recognition of mock threshold graphs is in Class P.

Something strange happens when we ask planarity questions about mock threshold graphs.  Outerplanarity and planarity become very different.

\begin{theorem}[Behr--Sivaraman--Zaslavsky \cite{BSZp}]\label{mtgp}
Outerplanar mock threshold graphs are characterized by finitely many forbidden topological induced subgraphs.  
Planar mock threshold graphs are characterized by infinitely many forbidden topological induced subgraphs.  
\end{theorem}

\section{Induction-Hereditary Classes}

\begin{question}[Sivaraman]
What can we say about a general hereditary class in the induced-subgraph ordering $\sqsubseteq$?
\end{question}

This is rather vague, and so general there might be nothing.  But there is something!  Given an induction-hereditary class $\cG$, define the sequence of cardinalities $\Phi(\cG) = (\phi_n(\cG))_{n=1}^\infty$ by $\phi_n(\cG) = |\{H \in \Forb(\cG) : |V(H)| = n \}|$.

\begin{theorem}[Sivaraman \cite{Siv}]
There exist positive real numbers $a, b$ such that every sequence $\{f_n\}$ with $f_n \leq a e^{bn}$ is the sequence $\Phi(\cG)$ for some induc\-tion-hereditary class $\cG$.
\end{theorem}


\end{document}